# KIOPS: A fast adaptive Krylov subspace solver for exponential integrators

Stéphane Gaudreault [a,*], Greg Rainwater [b], Mayya Tokman [b]

[a] *Recherche en prévision numérique atmosphérique, Environnement et Changement climatique Canada, 2121 Route Transcanadienne, Dorval, Québec, H9P 1J3, Canada*

[b] *School of Natural Sciences, University of California, 5200 N. Lake Road, Merced, CA 95343, United States*

ARTICLE INFO

ABSTRACT

This paper presents a new algorithm KIOPS for computing linear combinations of $\varphi$-functions that appear in exponential integrators. This algorithm is suitable for large-scale problems in computational physics where little or no information about the spectrum or norm of the Jacobian matrix is known *a priori*. We first show that such problems can be solved efficiently by computing a single exponential of a modified matrix. Then our approach is to compute an appropriate basis for the Krylov subspace using the incomplete orthogonalization procedure and project the matrix exponential on this subspace. We also present a novel adaptive procedure that significantly reduces the computational complexity of exponential integrators. Our numerical experiments demonstrate that KIOPS outperforms the current state-of-the-art adaptive Krylov algorithm *phipm*.

## 1. Introduction

Exponential integration received increased attention recently as an efficient alternative strategy to standard methods for solving systems of ordinary differential equations (ODE). Exponential integrators have the advantage of being accurate and, similarly to implicit methods, possess good stability properties, allowing integration with large time steps. These methods involve computation of an exponential or an exponential-like function of a Jacobian matrix or an approximation to it (e.g. see review articles [1,2]). Approximation of such matrix functions or their products with vectors constitutes the main computational cost of an exponential integrator. Typically it is the latter approximation, i.e. product of exponential-like matrix functions with vectors, that is required for an implementation of an exponential method.

A number of methods have been proposed to calculate the exponential or the exponential-like functions of a matrix or their product with a given vector (see [3] for a review). Most of them, however, are of little practical use for large-scale stiff matrices due either to high computational cost or numerical stability issues. These challenges are described in review papers of Moler and Van Loan [4,5]. Other algorithms that are more suitable for large stiff matrices [6,7] require some information about the norm or the spectrum of a matrix. It is common, however, that the matrix in question is not given in the explicit form and only matrix-vector products can be calculated. This is the case, for instance, when the matrix is a Jacobian that results from a complicated spatial discretization of a system of partial differential equations (PDE). As a result of these considerations, a conclusion can be drawn that the Krylov subspace projection-based techniques are among the most promising methods for problems where exponential or exponential-like functions of a large stiff matrix have to be computed and little or no information can be elicited about this matrix *a priori*.

* Corresponding author.
*E-mail addresses:* stephane.gaudreault2@canada.ca (S. Gaudreault), grainwater@ucmerced.edu (G. Rainwater), mtokman@ucmerced.edu (M. Tokman).

The basic idea of the Krylov subspace approach is to project the exponential of a large matrix onto a relatively small Krylov subspace where calculating the exponential is significantly less computationally expensive. Recent progress in computational linear algebra has led to efficient Krylov subspace algorithms such as the EXPOKIT software of Sidje [8], restarted Krylov methods [9–12], block Krylov subspaces [13–16], time-parallel methods [17,15], the shift-and-invert acceleration [18–20] and the adaptive methods [21,22]. The *phipm* adaptive method of Niesen and Wright [21] has been shown to be the most efficient option for the problems under consideration [23].

For the large scale geophysical fluid dynamics problems that motivate our study, it was found by Clancy and Pudykiewicz [24] that semi-implicit predictor–corrector schemes [25] where still more efficient. Gaudreault and Pudykiewicz [26] further analyzed these results and concluded that the computational cost of the Arnoldi procedure [27] was the origin of the issue. The performance was then improved by using an incomplete orthogonalization instead of the Arnoldi procedure based on full orthogonalization. This technique has been successfully used for the simulation of the shallow water equations on the sphere with second and third order exponential propagation iterative (EPI) schemes [28].

The incomplete orthogonalization procedure (hereafter denoted as IOP) was originally proposed by Saad as an eigenvalue algorithm for general non-symmetric matrices [29] and to solve systems of linear equations [30]. The application of IOP to approximating the matrix exponential was studied recently by Koskela [31] and by Vo and Sidje [32].

The aim of this article is to explore this technique further with a particular focus on the efficient calculation of $\varphi$-functions within exponential integrators. We call the resulting algorithm the Krylov with Incomplete Orthogonalization Procedure Solver (KIOPS). This new method has been carefully designed to allow for an efficient implementation of single or multi-stage exponential integrators, such as those recently proposed by Rainwater and Tokman [33].

The paper is organized as follows. In Section 2 we describe the main application of KIOPS method, i.e. evaluation of linear combinations of $\varphi$ functions within exponential integrators. Section 3.1 presents a theorem that shows a way to evaluate a linear combination of the $\varphi$-functions used in our KIOPS algorithm. The KIOPS algorithm is presented in Sections 3.2–3.7 where we describe how the method uses Krylov subspace projection with incomplete orthogonalization and a new adaptivity procedure. Our numerical experiments in Section 4 validate the performance of the proposed solver and demonstrate its relative efficiency compared to *phipm*. Conclusions are presented in Section 5.

## 2. Approximating $\varphi$-functions within exponential integrators

The main application of KIOPS and other algorithms that approximate products of $\varphi$-functions and vectors is their use within an exponential integrator designed to solve initial value problems for large scale stiff systems of ordinary differential equations (ODE) of the form

$$\begin{aligned} &\frac{d}{dt}u(t) = f(u(t)), \\ &u(t_0) = u_0, \qquad t \in [t_0, t_{end}], \quad u, f(u) \in \mathbb{R}^N. \end{aligned} \tag{1}$$

Differential equations of this form arise in many contexts in the natural and social sciences and engineering disciplines. In most applications, the independent variable $t$ usually represents time, $N$ is the number of degrees of freedom, the vector-valued function $u(t)$ represents some unknown dynamical quantities and $f$ is a vector-valued function describing all forces driving the system.

Exponential integrators gained much interest in recent years as an efficient alternative to implicit schemes in integrating stiff systems and a number of exponential methods have been introduced (e.g. [34–40,28,41], see also review [2]). When any exponential integrator is applied to a complex nonlinear problem its main computational cost is the evaluation of the exponential matrix-functions vectors products. The new KIOPS technique applies to any exponential integrator that employs approximation of these products. To make the computational savings of KIOPS concrete, however, in this paper we will focus on a class of exponential propagation iterative methods of Runge–Kutta-type (EPIRK). These schemes are particularly designed to gain efficiency by exploiting the properties of an algorithm used for the evaluation of $\varphi$-functions [23,33]. Specifically, EPIRK schemes proposed in [23,33] are particularly efficient when used with an adaptive Krylov method [42,43].

To illustrate what are the computational costs associated with approximating $\varphi$-functions within an exponential integrator we consider a three-stage EPIRK method that can be written as

$$\begin{aligned} U_{n2} &= u_n + \alpha_{11}\psi_{11}(g_{11}h_n J_n)h_n f(u_n) \\ U_{n3} &= u_n + \alpha_{21}\psi_{21}(g_{21}h_n J_n)h_n f(u_n) + \alpha_{22}\psi_{22}(g_{22}h_n J_n)h_n r(U_{n2}) \\ u_{n+1} &= u_n + \beta_1\psi_{31}(g_{31}h_n J_n)h_n f(u_n) + \beta_2\psi_{32}(g_{32}h_n J_n)h_n r(U_{n2}) \\ &\quad + \beta_3\psi_{33}(g_{33}h_n J_n)h_n(-2r(U_{n2}) + r(U_{n3})) \end{aligned} \tag{2}$$

where $h_n = t_{n+1} - t_n$, $r(u) = f(u) - f(u_n) - J_n(u - u_n)$ is the nonlinear remainder of the second order Taylor expansion of $f(u)$ and $\psi_{ij}(z)$ are linear combinations of exponential-like functions defined as

$$\psi_{ij}(z) = \sum_{k=1}^{K} p_{ijk}\varphi_k(z), \quad \varphi_k(z) = \int_0^1 e^{z(1-\theta)} \frac{\theta^{k-1}}{(k-1)!} d\theta. \tag{3}$$

Constant coefficients $\alpha_{ij}$, $\beta_{ij}$ and $p_{ijk}$ can be derived using either classical or stiff order conditions and methods of up to order five have been proposed in [44,33]. Note that since we are interested in problems where $J_n \in \mathbb{R}^{N\times N}$ with $N >> 1$, the largest per-time step computational cost of (2) lies in evaluating the matrix function-vector products $\psi_{ij}(gA)b$ ($g \in \mathbb{R}$, $A \in \mathbb{R}^{N\times N}$, $b \in \mathbb{R}^N$) for different vectors $b$ (e.g. $b = h_n f(u_n)$, $b = h_n r(U_{n2})$, etc.). Thus the most efficient methods should have coefficients $p_{ijk}$, $\alpha_{ij}$ and $\beta_j$ that both reduce the number of $\psi_{ij}(gA)b$ terms required as well as make each of these evaluations as computationally cheap as possible. A strategy to achieve this goal was proposed in previous publications [23,33] and involves exploiting the structure of the adaptive Krylov projection algorithm *phipm* [21].

The *phipm* algorithm evaluates linear combinations

$$w(\tau) = \varphi_0(\tau A)b_0 + \tau\varphi_1(\tau A)b_1 + \tau^2\varphi_2(\tau A)b_2 + \ldots + \tau^p\varphi_p(\tau A)b_p, \tag{4}$$

where $\tau \in \mathbb{R}$, $A \in \mathbb{R}^{N\times N}$ and $b_0, b_1, \ldots, b_p \in \mathbb{R}^N$. The *phipm* algorithm proceeds by considering (4) as a solution to the initial value problem

$$w'(\tau) = Aw(\tau) + b_1 + \tau b_2 + \ldots + \frac{\tau^{p-1}}{(p-1)!} b_p \quad w(0) = b_0 \tag{5}$$

The algorithm evaluates solution to (5) over $k$ subintervals $0 \le \tau_1 \le (\tau_1 + \tau_2) \le \ldots \le \sum_{l=1}^{k} \tau_l = \tau$ using a substepping procedure. Each substep $i$ is comprised of two main parts. The first part consists of computing the approximation

$$w(t_{i+1}) = \tau_i^p \varphi_p(\tau_i A)\tilde{w}_p + \sum_{j=0}^{p-1} \frac{\tau_i^j}{j!} \tilde{w}_j, \tag{6}$$

where $t_{i+1} = \sum_{l=1}^{i} \tau_l$ and the vectors $\tilde{w}_j$ are calculated recursively from the relation

$$\tilde{w}_0 = w(t_i) \quad \text{and} \quad \tilde{w}_j = A\tilde{w}_{j-1} + \sum_{\ell=0}^{p-j} \frac{\tau_i^\ell}{\ell!} b_{j+\ell}, \quad j = 1, \ldots, p. \tag{7}$$

The second part of each substep is the adaptive algorithm that selects the substep size $\tau_{i+1}$ and the dimension $m_{i+1}$ of the Krylov projection to be used for the next substep. To approximate

$$\varphi_0(A)b_0 + \varphi_1(A)b_1 + \varphi_2(A)b_2 + \ldots + \varphi_p(A)b_p, \tag{8}$$

the substepping procedure should be performed over the whole interval $\tau \in [0, 1]$.

Note that if only a single $\varphi_q(\tau A)b_q$ is involved in the linear combination (4) (i.e. $b_j = 0$ for $j \ne q$) then the substepping procedure can be modified to compute several terms of type $\varphi_q(T_i A)b_q$ for several values of $T_i$ simultaneously. This is accomplished by ensuring that the stops at each $\tau = T_i$ are made and the result is recorded and scaled by $1/T_i^q$.

Efficient EPIRK method (2) is derived when the terms $\psi_{ij}(gA)v$ can be combined into a minimum number of groups so that all such terms in each group are evaluated using only one execution of the adaptive Krylov projection algorithm. For example, a stiffly accurate fourth-order scheme EPIRK4s3 [45,43]

$$\begin{aligned}
U_{n2} &= u_n + \frac{1}{8}\varphi_1\left(\frac{1}{8}h_n J_n\right) h_n f(u_n), \\
U_{n3} &= u_n + \frac{1}{9}\varphi_1\left(\frac{1}{9}h_n J_n\right) h_n f(u_n), \\
u_{n+1} &= u_n + \varphi_1\left(h_n J_n\right) h_n f(u_n) \\
&\quad + \left(1892\varphi_3\big(h_n J_n\big) - 42336\varphi_4\big(h_n J_n\big)\right) h_n r(U_{n2}) \\
&\quad + \left(1458\varphi_3\big(h_n J_n\big) - 34992\varphi_4\big(h_n J_n\big)\right) \\
&\qquad \cdot h_n \left(r(U_{n3}) - 2r(U_{n2})\right)
\end{aligned} \tag{9}$$

requires only two calls to *phipm*. One call is used to evaluate both terms $\varphi_1\left(\frac{1}{8}h_n J_n\right) h_n f(u_n)$ and $\varphi_1\left(\frac{1}{9}h_n J_n\right) h_n f(u_n)$. Note that for this approximation we need to only substep solution (4) over the interval $\tau \in [0, \frac{1}{8}]$. The *phipm* approximation allows us to compute the linear combination (4) which involves coefficients $\tau^i$. The terms we need in (9) do not have these coefficients. Since the terms in (9) involve only single $\varphi_1$ term with all other $b_i = 0$ ($i \ne 1$), this problem can be easily remedied by multiplying the results of *phipm* by factors $1/\tau^i$. Such fix would not work if several $\varphi_i$'s were involved unless we need the final value of $\tau$ be 1. This is precisely the case for the second call to *phipm* to evaluate the linear combination in the last stage of (9).

The second and last call to *phipm* is used to evaluate the linear combination

$$\varphi_0(\tau h_n J_n)b_0 + \tau\varphi_1(\tau h_n J_n)b_1 + \tau^2\varphi_2(\tau h_n J_n)b_2 + \tau^3\varphi_3(\tau h_n J_n)b_3 + \tau^4\varphi_4(\tau h_n J_n)b_4 \tag{10}$$

with

$$\begin{aligned} b_1 &= h_n f(u_n), \\ b_2 &= 0, \\ b_3 &= 1892 h_n r(U_{n2}) + 1458 h_n\,(r(U_{n3}) - 2r(U_{n2})) \\ b_4 &= -42336 h_n r(U_{n2}) - 34992 h_n\,(r(U_{n3}) - 2r(U_{n2})) \end{aligned} \tag{11}$$

for $\tau = 1$. This linear combination involves all exponential terms involved in approximating $u_{n+1}$. Clearly the factors $1/\tau^i$ do not pose a difficulty in evaluation of this linear combination since we substep *phipm* over $\tau \in [0, 1]$.

Thus the structure of this exponential integrator utilizes *phipm* for two types of tasks:

(I) Computing several terms of type $\varphi_k(\tau A)b_k$ involving a single function $\varphi_k$ and a vector $b_k$ for several values of $\tau$;
(II) Approximating linear combinations (4) with $\tau = 1$.

Both higher-order higher-stage-number EPIRK schemes as well as other exponential integrators like exponential Rosenbrock methods [41] utilize similar structural properties to construct and implement exponential time integrators efficiently.

We have modified the original version of the *phipm* to enable these two options of computing either a single $\varphi_q(\tau A)b_q$ or the whole linear combination (4) for multiple values of $\tau = T_i$. In the modified *phipm*, we use the Algorithm 3 described in section 3.3.1 to keep track of intermediate solutions and we scale the final results by $1/\tau^i$ whenever task I is executed. The rest of the procedure is identical to the one described in [21]. In the subsequent text every mention of *phipm* refers to this slightly modified version of the original *phipm* algorithm.

For general problems where matrix $A$ is not explicitly available and little is known about its spectrum, *phipm* has been demonstrated to be the most efficient algorithm to date for the implementation of exponential integrators [23]. This algorithm has, however, significant drawbacks and the following three considerations were the main driving force behind our search of alternative techniques.

- The convergence of the *phipm* algorithm is often inconsistent. For example, it may happen that computing a solution with a small error is significatively faster than computing a solution with a lower accuracy. We also noticed that this behavior may occur alternately when more and more accurate solutions are being calculated. As we will see later in the section presenting our numerical experiments, this results in zig-zags in a precision diagram.
- The substepping of equation (5) requires several explicit multiplication by the matrix $A$ in Eq. (7). As discussed in [21], [22] and [26], this is not only costly, but the procedure might also become increasingly sensitive to rounding errors as $p$ increases. Since the last stage of a high-order EPIRK scheme often involve multiplications by relatively large coefficients, this sensitivity to rounding errors could be worrying if $\|A\|$ is large.
- The adaptive procedure of *phipm* is based on somewhat rough estimates of the floating point operations count on a single processor machine. Important details related to a specific computer architecture are not taken into account.

## 3. The KIOPS algorithm

The KIOPS algorithm, outlined in Algorithm 1, builds on the ideas of *phipm*, but modifies both the substepping procedure and the adaptive algorithm to significantly improve efficiency and accuracy of approximating (4). Specifically:

- Instead of substepping equation (5) we extend theoretical result of Sidjie [8] to linear combinations of $\varphi$-functions (4) and use exponential of an augmented matrix to compute all terms $\tau^i\varphi_i(\tau A)b_i$ simultaneously using one Krylov projection.
- We use incomplete rather than full orthogonalization procedure for Krylov projections.
- We propose a different adaptivity method to select $\tau_i$ and Krylov subspace sizes $m_i$ which brings more efficiency to the overall approximation of (4).

Below we provide a detailed description of the KIOPS algorithm and demonstrate its efficiency and accuracy within exponential integrators compared to *phipm* using a set of numerical examples.

### 3.1. Computing linear combinations of $\varphi$-functions

As mentioned before, the *phipm* requires several explicit multiplications by $A$ in Eq. (7). Hence, it becomes increasingly sensitive to round off errors as $p$ increases. The alternative approach presented in this section allows to evaluate all terms of the linear combination (4) simultaneously and thus completely replaces the substepping of (5).

Once again consider the task of approximating

$$\varphi_0(\tau A)b_0 + \tau\varphi_1(\tau A)b_1 + \tau^2\varphi_2(\tau A)b_2 + \ldots + \tau^p\varphi_p(\tau A)b_p \tag{12}$$

where $\tau \in \mathbb{R}$, $A \in \mathbb{R}^{N\times N}$ and $b_0, b_1, \ldots, b_p \in \mathbb{R}^N$. In a typical application, the matrix $A$ is large and sparse. Exact evaluation of individual $\varphi$-functions and vector products is then prohibitively computationally expensive. We choose a more efficient method to compute the linear combination (12) using a single exponential of a slightly larger matrix [8,6,22]. The following theorem is an extension of the result obtained in [8] for the case $b_i = c$ for $i = 0, \ldots, p$. Theorem 1 shows how the problem of computing (12) can be solved by computing a single exponential of a matrix. Hereafter, we denote by $I_l$ the $l \times l$ identity matrix, $e_p = (0, \ldots, 0, 1)^\intercal \in \mathbb{R}^p$ is the last canonical basis vector in $\mathbb{R}^p$. The colon operator "$a : b$" is an operation that generates the indices ranging from $a$ to $b$. This operator is used to represent a subset of the elements from a vector or a matrix and it has the lowest priority.

**Theorem 1.** *Let $A \in \mathbb{R}^{N\times N}$, $B = [b_p, \ldots, b_2, b_1] \in \mathbb{R}^{N\times p}$, $K = \begin{bmatrix} 0 & I_{p-1} \\ 0 & 0 \end{bmatrix} \in \mathbb{R}^{p\times p}$, $v = [b_0^\intercal, e_p^\intercal]^\intercal \in \mathbb{R}^{N+p}$ and $\tau \in \mathbb{R}$. We define the augmented matrix*

$$\tilde{A} = \begin{bmatrix} A & B \\ 0 & K \end{bmatrix} \in \mathbb{R}^{(N+p)\times(N+p)} \tag{13}$$

*and let $w = e^{\tau\tilde{A}}v$. Then, the first $N$ elements of the vector $w$ are given by*

$$w(1:N) = \sum_{j=0}^{p} \tau^j \varphi_j(\tau A) b_j \tag{14}$$

*and the rest of the vector is*

$$w(N+1:N+p) = \left[\frac{\tau^{p-1}}{(p-1)!}, \ldots, \tau, 1\right]^\intercal \tag{15}$$

**Proof.** Since $\tilde{A}$ is block upper triangular, its exponential has the form

$$e^{\tau\tilde{A}} = \begin{bmatrix} e^{\tau A} & E_{12} \\ 0 & e^{\tau K} \end{bmatrix}. \tag{16}$$

From [8], we have

$$e^{\tau K} = \begin{bmatrix} 1 & \frac{\tau}{1!} & \cdots & \frac{\tau^{p-2}}{(p-2)!} & \frac{\tau^{p-1}}{(p-1)!} \\ 0 & 1 & \ddots & & \frac{\tau^{p-2}}{(p-2)!} \\ \vdots & & \ddots & \ddots & \vdots \\ & & & 1 & \frac{\tau}{1!} \\ 0 & & \cdots & 0 & 1 \end{bmatrix} \tag{17}$$

and the columns of $E_{12}$ are given by the Theorem 2.1 of Al-Mohy and Higham [6], with $\ell = 0$, as

$$E_{12}(1:N, N+j) = \sum_{k=1}^{j} \tau^k \varphi_k(\tau A) b_k\,, \quad \text{for } j = 1 \text{ to } p. \tag{18}$$

Inserting Eqs. (17) and (18) into Eq. (16), we obtain the following matrix

$$e^{\tau \tilde{A}} = \begin{bmatrix} e^{\tau A} & \tau\varphi_1(\tau A)b_1 & \tau\varphi_1(\tau A)b_1 + \tau^2\varphi_2(\tau A)b_2 & \dots & \sum_{k=1}^{p} \tau^k \varphi_k(\tau A) b_k \\ 0 & 1 & \frac{\tau}{1!} & \dots & \frac{\tau^{p-1}}{(p-1)!} \\ \vdots & & \ddots & \ddots & \vdots \\ & & & 1 & \frac{\tau}{1!} \\ 0 & & & 0 & 1 \end{bmatrix}. \tag{19}$$

It only remains to multiply (19) by $v = \left[b_0, e_p\right]^\intercal$ to obtain the vector $w$. □

Now we will use this result inside an iterative procedure to evaluate $e^{\tau \tilde{A}}$.

### 3.2. Krylov adaptive method

The KIOPS algorithm is structured around the idea of utilizing Theorem 1 to approximate

$$w(\tau) = \sum_{j=0}^{p} \tau^j \varphi_j(\tau A) b_j \tag{20}$$

using the exponential of an augmented matrix. This matrix exponential can be computed using a polynomial approximation of the form

$$w(\tau) = e^{\tau \tilde{A}} v \approx P_{m-1}(\tau \tilde{A}) v, \tag{21}$$

where $P_{m-1}$ is a polynomial of degree $m-1$. There are a number of methods that employ a polynomial approximation of the form of Eq. (21) including truncated Taylor series approximation, Leja interpolation, and Chebyshev polynomials-based algorithms. The disadvantage of most of these methods is that they require information about the spectrum or norm of the matrix. As mentioned earlier obtaining this extra information about the matrix can be impossible or prohibitively expensive computationally especially for problems where the action of the matrix-vector product is given by an external "matvec" subroutine. In this work, we consider only approaches that do not require any knowledge about the norm or spectrum of the matrix. We also avoid methods that require inversion of $\tilde{A}$ because this matrix is singular.

Since the approximation (21) is an element of the $m$-dimensional Krylov subspace

$$K_m(\tilde{A}, v) = \text{span}\{v, \tilde{A}v, \dots, \tilde{A}^{m-1}v\}, \tag{22}$$

the problem can be recast into the search for an element of $K_m$ that approximates $w(\tau)$ [46]. The approximation of the vector $w(\tau)$ by an element of a Krylov subspace is made up of two important steps. The first step is the computation of an appropriate basis for the Krylov subspace and of a smaller matrix that represents the projection of the action of $\tilde{A}$ on this subspace. Then, in the second step, the matrix exponential of the smaller matrix is computed using a standard technique and the result is projected back onto the original large space. Notice that the definition (22) does not include the factor $\tau$ since the Krylov subspace associated with $\tilde{A}$ and $\tau\tilde{A}$ are the same for any $\tau \in \mathbb{R}$.

We will see in section 3.4 that if we want to obtain good accuracy, then the size of the Krylov space $m$ has to be large when $\|\tau\tilde{A}\|$ is large. This is worrying because it may indicate that an impractical amount of memory storage and computational cost could be necessary to obtain a sufficiently small error. A more efficient approach is to apply the Krylov subspace method iteratively as in the work of Sidje [8]. The key idea is that computing the action of the matrix exponential is equivalent to solving a linear ODE to split $\tau$ into a sum of smaller intervals, such that

$$e^{\tau \tilde{A}} v = e^{(\tau_1 + \tau_2 + \dots + \tau_k)\tilde{A}} v \tag{23}$$

$$= e^{\tau_k \tilde{A}}(\dots(e^{\tau_2 \tilde{A}}(e^{\tau_1 \tilde{A}} v))). \tag{24}$$

Thus, the problem of computing the matrix exponential

$$w(\tau) = e^{\tau \tilde{A}} v \tag{25}$$

is equivalent to finding a solution for the initial value problem

$$\frac{d}{d\tau} w(\tau) = \tilde{A} w(\tau), \qquad w(0) = v. \tag{26}$$

From this viewpoint, one could iteratively solve the recurrence

$$w(0) = v \tag{27}$$

$$w(\tau_{i+1}) = e^{\tau_i \tilde{A}} w(\tau_i) \tag{28}$$

and use an adaptive procedure in order to control the errors of the method.

Given a function that computes the product $Av$ of a matrix $A \in \mathbb{R}^{N \times N}$ and a vector $v \in \mathbb{R}^N$ each substep $\tau_l$ the KIOPS algorithm proceeds according to the following steps:

(i) Given current substep size $\tau_l$ and the size of the Krylov subspace $m_l$ an incomplete orthogonalization-based Krylov projection algorithm is executed to compute $e^{\tau_l \tilde{A}}$ defined in (19) (this step is described in detail in sections 3.3 and 3.3.1).

(ii) The local error estimate is computed and the approximation from step (i) is accepted for the current interval only if the error is within user-specified tolerance (see sections 3.4, 3.5).

(iii) If the approximation in step (ii) is accepted, the adaptive algorithm then determines whether it is more cost efficient to compute the approximation at the next substep with a larger value of $\tau_{l+1}$ or a smaller Krylov subspace size $m_{l+1}$. If the approximation in step (ii) was rejected, the adaptive algorithm decides whether it is more cost efficient to obtain a better estimate by reducing the substep size $\tau_l$ or increasing the Krylov subspace size $m_l$ (see sections 3.6.1 and 3.6.2).

As seen in the previous section, an exponential integrator typically requires either evaluation of $w$ at several values of $\tau = g_{ij} = T_i < 1$ (Task I, usually in cases of $w$ composed of a single $\varphi_j$-function) or over the interval $[0, \tau = 1]$ (Task II). When several values of $\tau$ are needed, they are provided as an argument to Algorithm 1 in the array $T = [T_1, \ldots, T_{\text{end}}]$. Algorithm 1 then returns several linear combinations of the form

$$\varphi_0(T_1 A)b_0 + \varphi_1(T_1 A)b_1 + \varphi_2(T_1 A)b_2 + \ldots + \varphi_p(T_1 A)b_p$$

$$\ldots$$

$$\varphi_0(T_{\text{end}} A)b_0 + \varphi_1(T_{\text{end}} A)b_1 + \varphi_2(T_{\text{end}} A)b_2 + \ldots + \varphi_p(T_{\text{end}} A)b_p.$$

**Algorithm 1** KIOPS: Evaluate linear combination (12).

1: input: $T = [T_1, \ldots, T_{\text{end}}]$, $A \in \mathbb{R}^{N \times N}$, $U = \left[b_p, \ldots, b_2, b_1, b_0\right]$, tol (default $1e-7$), $m_{\text{init}}$ (default 10), $m_{\text{min}}$ (default 10), $m_{\text{max}}$ (default 128), Task (I or II)
2: $\tau = T_{\text{end}}$,
3: $m = \max(m_{\text{min}}, \min(m_{\text{init}}, m_{\text{max}}))$
4: $T_{\text{now}} = 0$
5: $j = 0$
6: $\ell = 0$
7: $n_{\text{steps}} = \text{length}(T)$
8: $w(1:N, 1) = b_0$
9: **while** $T_{\text{now}} < T_{\text{end}}$ **do**
10: **if** j == 0 **then**
11: {Compute the first Krylov basis vector}
12: $H = 0$
13: $w(N+1, N+p, \ell) = \left[\frac{(T_{\text{now}})^{p-1}}{(p-1)!}, \ldots, T_{\text{now}}, 1\right]^{\mathsf{T}}$
14: $V(:,1) = \frac{1}{\|w(:,\ell)\|} w(:,\ell)$
15: **end if**
16: V, H, j = Algorithm 2 ($A$, $\left[b_p, \ldots, b_2, b_1\right]$, $V$, $j$, $m$) {section 3.3}
17: $H(1, j+1) = 1$
18: $F = \exp(\tau H(1:j+1, 1:j+1))$ {section 3.3.1}
19: Compute the local error estimate as in section 3.4.
20: Calculate suggested $\tau_{\text{new}}$ and $m_{\text{new}}$ as in section 3.6.
21: Choose to change $m$ for $m_{\text{new}}$ or $\tau$ for $\tau_{\text{new}}$.
22: **if** Acceptance criterion described in section 3.5 is satisfied **then**
23: $w$, $\ell$ = Algorithm 3 ($T$, $H$, $w$, $\ell$) { Update the $w$ vector }
24: $T_{\text{now}} = T_{\text{now}} + \tau$
25: $j = 0$
26: **else**
27: $H(1, j+1) = 0$ { Restore the original matrix }
28: **continue**
29: **end if**
30: **end while**
31: **if** TASK I **then**
32: **for** $l = 1$ **to** $n_{\text{steps}}$ **do**
33: $w(:,l) = w(:,l) * (1/T(l))^p$
34: **end for**
35: **end if**
36: **return** $w_\ell(1:N) = e^{T_\ell \tilde{A}} v$, for all $1 \leq \ell \leq \text{end}$, and $m$.

To simplify the notation in the subsequent sections we will refer to the substep size $\tau_l$ as simply $\tau$ and the Krylov substep size $m_l$ as $m$. Keep in mind, however, that these values change over each substep $l$.

### 3.3. Building a basis for the Krylov subspace

The first part of the Krylov subspace method consists of computing a set of basis vectors $\{v_1, \ldots, v_m\}$ for the $m$-dimensional Krylov subspace $K_m$. This is equivalent to finding a matrix

$$V_m = [v_1, \ldots, v_m] \tag{29}$$

whose column space is $K_m$. In a typical implementation, the Arnoldi procedure is used to compute this matrix. As discussed in [26], performing the Arnoldi procedure takes $O(m^2 \cdot [N+p])$ operations and constitutes the primary computational cost of the Krylov subspace projection technique. For this reason we turn to the incomplete orthogonalization procedure of length 2 whose time complexity is $O(m \cdot [N+p])$.

Starting with the vector $v_1 = v/\|v\|$, the incomplete orthogonalization procedure (Algorithm 2) produces the factorization

$$\tilde{A}\, V_m = V_m\, H_m + h_{m+1,m}\, v_{m+1}\, e_m^{\mathsf{T}}, \tag{30}$$

where $e_m = (0, \ldots 0, 1)^{\mathsf{T}} \in \mathbb{R}^m$ denotes the last canonical basis vector in $\mathbb{R}^m$. The most important by-products of this factorization are the matrix $V_m \in \mathbb{R}^{(N+p)\times m}$ and $H_m \in \mathbb{R}^{m\times m}$. The matrix $H_m$ has a banded structure and it can be seen as an oblique projection of the action of the matrix $\tilde{A}$ on the Krylov subspace. The entry $h_{m+1,m} \in \mathbb{R}$ can be interpreted as a kind of residual of the projection onto the Krylov subspace and it will enter into the formulation of an error estimate in section 3.4. Although we want to generate a basis of an $m$-dimensional space, the algorithm produces $m+1$ vectors in general. The last vector $v_{m+1} \in \mathbb{R}^{(N+p)}$ will not be used in our approximation scheme. There exists alternative corrected scheme [46] that makes use of this vector, but our experiments showed that it has slightly slower convergence than the standard scheme when used in this solver.

The main difference between Algorithm 2 and an equivalent algorithm based on Arnoldi procedure is that each new vector is orthogonalized only against the two previous ones instead of all of them. Hence, the matrix $V_m$ obtained after Algorithm 2 has rank $m$ and its columns span $K_m$, but they do not form an orthonormal basis in general.

**Algorithm 2** Incomplete orthogonalization procedure of length 2.

1: **Input:** $A \in \mathbb{R}^{N\times N}$, $B \in \mathbb{R}^{N\times p}$, $V \in \mathbb{R}^{N+p\times m_{\max}+1}$, $j$, $m$
2: **while** $j < m$ **do**
3: $\quad j = j+1$
4: $\quad V(1:N, j+1) = A \cdot V(1:N, j) + B \cdot V(N+1:N+p, j)$
5: $\quad V(N+1:N+p-1, j+1) = V(N+2:N+p, j)$
6: $\quad V(N+p, j+1) = 0$
7: $\quad$ **for** $i = \max(1, j-1)$ **to** $j$ **do**
8: $\qquad H(i, j) = V(:, i)^{\mathsf{T}} \cdot V(:, j+1)$
9: $\qquad V(:, j+1) = V(:, j+1) - H(i, j) \cdot V(:, i)$
10: $\quad$ **end for**
11: $\quad s = \|V(:, j+1)\|$
12: $\quad$ **if** $s \approx 0$ **then**
13: $\qquad$ happy_breakdown = true
14: $\qquad$ **break**
15: $\quad$ **end if**
16: $\quad H(i+1, j) = s$
17: $\quad V(:, j+1) = \frac{1}{s} V(:, j+1)$
18: **end while**
19: **return** $V$, $H$, $j$

It is worth noting that Algorithm 2 does not require the assembly of the augmented matrix $\tilde{A}$. In effect, the product $\tilde{A} v_j$ is formulated, at the lines 4 to 6, using block multiplication involving only the action of $A$ and the matrix $B$. This formulation allows the action of the matrix $A$, or some approximation of it, to be provided by an external subroutine. Typically, such subroutine tries to exploit the sparsity patterns of $A$ to reduce the computational cost associated with the matrix-vector product.

#### 3.3.1. Approximation of the exponential in the Krylov subspace

Since columns of $V_m$ obtained using IOP do not form an orthonormal basis, $V_m^T V_m$ is not an identity matrix as is the case for the Arnoldi procedure. We can nevertheless use arguments similar to those presented by Saad [46] and by Higham [3] to obtain a theorem that will justify our approximation scheme.

**Theorem 2.** *Let $\tilde{A} \in \mathbb{R}^{(N+p)\times(N+p)}$, $v \in \mathbb{R}^{(N+p)}$ and let $V_m \in \mathbb{R}^{(N+p)\times m}$, $H_m \in \mathbb{R}^{m\times m}$ be the result of m steps of Algorithm 2. Then for any polynomial $P_j$ of degree $j < m-1$ we have*

$$P_j(\tilde{A})v = \|v\|\, V_m\, P_j(H_m)\, e_1 \tag{31}$$

**Proof.** Without loss of generality, the proof can be done by induction on $j$ for the polynomials of the form $P_j(\tilde{A}) = \tilde{A}^j$. For $j=0$, we have that $\tilde{A}^0 v = v = \|v\|\, V_m\, H_m^0\, e_1$ since $V_m e_1 = v_1$. Assume that (31) is true for all $j \le m-2$, then

$$\begin{aligned} \tilde{A}^{j+1} v &= \tilde{A} \cdot \tilde{A}^j v = \tilde{A} V_m H_m^j (\|v\|\, e_1) \\ &= \left(V_m H_m + h_{m+1,m} v_{m+1}\, e_m^{\mathsf{T}}\right) H_m^j (\|v\|\, e_1) \\ &= \|v\|\, V_m H_m^{j+1}\, e_1 + h_{m+1,m} v_{m+1}\, e_m^{\mathsf{T}} H_m^j (\|v\|\, e_1) \end{aligned}$$

and using the fact that $h_{m+1,m} v_{m+1} e_m^{\mathsf{T}} H_m^j (\|v\|\, e_1) = 0$ whenever $j \le m-2$, we obtain that $\tilde{A}^{j+1} v = \|v\|\, V_m H_m^{j+1}\, e_1$, as required. □

This theorem, combined with the Taylor series definition of the exponential

$$e^{\tau \tilde{A}} = \sum_{k=0}^{\infty} \frac{1}{k!} (\tau \tilde{A})^k, \tag{32}$$

suggests the use of the following approximation scheme for the matrix exponential:

$$e^{\tau \tilde{A}} v \approx \|v\|\, V_m\, e^{\tau H_m} e_1. \tag{33}$$

Since in general $m \ll (N+p)$, the matrix exponential $e^{\tau H_m}$ can be computed using any standard method with dense output [4,5]. In this work, we use a diagonal Padé approximation combined with a scaling and squaring algorithm [47].

Algorithm 3 presents the procedure to calculate the solution for multiple stepsizes.

**Algorithm 3** Update the solution using Eq. (33).

1: input: $T = [T_1, \ldots, T_{\text{end}}]$, $H$, $w$, $\ell$
2: $n_\tau = 0$
3: $T_{\text{next}} = T_{\text{now}} + \tau$
4: **if** TASK I **then**
5:   **for** $k = \ell$ **to** $n_{\text{steps}}$ **do**
6:     **if** $\|T(k)\| < \|T_{\text{next}}\|$ **then**
7:       $n_\tau = n_\tau + 1$
8:     **end if**
9:   **end for**
10:   **if** $n_\tau > 0$ **then**
11:     $w(:, \ell + n_\tau) = w(:, \ell)$ { Copy current $w$ to $w$ we continue with. }
12:     **for** $k = 0$ **to** $n_\tau - 1$ **do**
13:       $\tilde{\tau} = T(\ell + k) - T_{\text{now}}$
14:       $F2 = \exp(\tilde{\tau} H(1:j, 1:j))$
15:       $w(:, \ell + k) = \|w(:, \ell)\| * V(1:N, 1:j) * F2(1:j, 1)$ { Using Eq. (33) }
16:     **end for**
17:     $\ell = \ell + n_\tau$ { Advance $\ell$. }
18:   **end if**
19: **end if**
    $w(:, \ell) = \|w(:, \ell)\| * V(1:N, 1:j) * F(1:j, 1)$ { Using Eq. (33) }
20: **return** $w$, $\ell$

### *3.4. Error estimates*

The approximation error for the exponential of a negative semidefinite matrix is bounded as in the following equation [31] (see also [32] for a comprehensive analysis of the convergence of IOP)

$$\|e^{\tau \tilde{A}} v - \beta V_m e^{H_m} e_1\| \le \frac{e^{\alpha(\tau \tilde{A})} \|\tau \tilde{A}\|^m + e^{\alpha(H_m)} \|H_m\|^m}{m!} \beta, \tag{34}$$

where $\alpha(\tilde{A}) = \max\{\mathrm{Re}(\lambda_i)\}$ denotes the spectral abscissa of $\tilde{A}$ (i.e. the supremum among the real part of the eigenvalues of $\tilde{A}$), $m$ is the size of the Krylov subspace, $\beta = \|v\|$ and $\|\tilde{A}\|$ is the matrix 2-norm. This *a priori* estimate is useful to gain insight on the convergence of the approximation scheme, but it cannot be used in practical calculations because it requires knowledge about the spectrum of $\tau \tilde{A}$.

A more practical *a posteriori* estimate can be derived analogously to the Theorem 5.1 of Saad [46]. This error estimate will be useful to formulate an acceptance criterion in section 3.5. The basic idea is that the error of the scheme (33) satisfies the expansion

$$e^{\tau \tilde{A}}\, v - \beta\, V_m\, e^{\tau H_m} e_1 = \beta\, \tau\, h_{m+1,m} \sum_{j=1}^{\infty} e_m^{\mathsf{T}}\, \varphi_j(\tau H_m) e_1\, (\tau \tilde{A})^{j-1}\, v_{m+1}. \tag{35}$$

If we assume that the magnitude of the terms of the series decreases rapidly as the solution converges, then the absolute value of the first term

$$\epsilon_m = \| \tau\, h_{m+1,m}\, e_m^{\mathsf{T}}\, \varphi_1(\tau H_m) \beta\, e_1 \| \tag{36}$$

can be used as an estimate of the error. This estimate can be computed with little extra cost. Invoking Theorem 1 again, the $\tau e_m^{\mathsf{T}} \varphi_1(H_m) e_1$ term can be obtained along with $e^{\tau H_m}$ in a single matrix exponential. To do that, we define the following $(m+1) \times (m+1)$ matrix

$$\tilde{H}_m = \begin{bmatrix} H_m & e_1 \\ 0 & 0 \end{bmatrix} \tag{37}$$

and evaluate its exponential

$$e^{\tau \tilde{H}_m} = \begin{bmatrix} e^{\tau H_m} & \tau \varphi_1(\tau H_m) e_1 \\ 0 & 1 \end{bmatrix}. \tag{38}$$

The term $\tau e_m^{\mathsf{T}} \varphi_1(H_m) e_1$ corresponds to the entry $(m, m+1)$ of the resulting matrix. It remains only to multiply by $\beta\, h_{m+1,m}$ to obtain the estimate (36).

Although no serious problems were found in this study, we note that this error estimate is less accurate when $m$ is large (see also the discussions [46] and [8]). Using a sharper error estimate or an appropriate residual notion (see e.g. [8,12,48]) could lead to better performance and accuracy. We intend to analyze this in detail in future studies.

### 3.5. Acceptance criterion

For an iteration to be accepted, it must satisfy a user defined tolerance. To verify this, we compute the scaled error

$$\omega = \frac{T_{\text{end}} \epsilon_m}{\tau\ \text{Tol}} \tag{39}$$

as in Niesen and Wright [21]. The step is accepted if $\omega <= \delta$, where $\delta = 1.4$ is a criterion intended to reduce the risk of rejection of the step. Since the adaptive procedure described in section 3.6.1 seeks to get a solution with a value of $\omega$ less than 0.9, there is little to no risk to allow the tolerance to be slightly exceeded occasionally. An iteration is rejected if $\omega > \delta$. The adaptive procedure can then extend the existing Krylov subspace basis with more vectors or use a smaller stepsize to improve the current solution. This flexibility to improve a missed iteration is an important difference with Krylov solvers for linear systems and eigenvalue problems where a bad iteration can wreak havoc.

### 3.6. Selection of parameters

The accuracy and the efficiency of our algorithm depend on the two important parameters: the size of the Krylov space $m_i$ and the substep size $\tau_i$. Thus their values must be chosen with care.

To start the algorithm, we use a somewhat optimistic first guess of $\tau = T_{\text{end}}$ and a user provided estimate $m = m_{\text{init}}$. A default value of $m_{\text{init}} = 10$ is used if no other value has been specified. However, experiments have indicated that a good estimate of the Krylov subspace size can drastically improve the performance of the method [26]. When the solver is used repetitively to integrate an ODE over many time steps, a possible strategy is to use the final value of $m_l$ from the previous time step as a first guess for the next one. This heuristic is justified when the nature of the problem to be solved in a particular time step resembles closely that from the previous step. Hence, the algorithm also returns the size of the Krylov subspace $m_l$ used in the last iteration as an output parameter. A solution is then produced with these parameters. If the acceptance criterion is satisfied, then our adaptive procedure will determine if $m$ or $\tau$ could be enlarged safely for the next substep. If, on the contrary, the acceptance criterion is not satisfied, then we must choose different parameters $m$ or $\tau$ and try again. The challenge is then to find optimal parameters to avoid step rejection and for quick convergence.

After each iteration (accepted or rejected), the adaptive procedure suggests a stepsize $\tau_{\text{new}}$ and a dimension of the Krylov subspace $m_{\text{new}}$ as in Algorithm 4. Two scenarios can then be considered and the most economical is selected:

(A) keeping $m$ constant and changing $\tau$ to $\tau_{\text{new}}$,
(B) keeping $\tau$ constant and changing $m$ to $m_{\text{new}}$.

Starting at a given substep and assuming that the new choice of parameter leads to acceptance for all of the subsequent steps we can see that the cost to complete the remaining steps of Algorithm 1 is bounded below by

$$C(\tau, m) = \left\lceil \frac{t_{\text{end}} - t_{\text{now}}}{\tau} \right\rceil \left( \text{cost}_{\text{iop}}(m) + \text{cost}_{\text{exp}}(N_{mult}, m) \right). \tag{40}$$

In (40) $\text{cost}_{\text{iop}}(m)$ is the cost of making $m$ steps of the Algorithm 2 (the incomplete orthogonalization procedure) that can be estimated as

$$\text{cost}_{\text{iop}}(m) = m[(2 \cdot \text{nz}(\tilde{A}) + (2p-1)N + p) + 8(N+p)] - 4(N+p) + 3(N+p)m \tag{41}$$

and $\text{cost}_{\text{exp}}(N_{mult}, m)$ is the cost of computing the dense matrix exponential given in [47] as

$$\text{cost}_{\text{exp}}(N_{mult}, m) = N_{mult}[2(m+p) - 1](m+p)^2, \tag{42}$$

where $\text{nz}(\tilde{A})$ denotes the number of nonzero elements in $\tilde{A}$ and $N_{mult}$ is the number of multiplications required for the scaling and squaring and Padé algorithms [47] that was selected for this study.

We notice that the cost of completing the algorithm with scenario (B), $C(\tau, m_{\text{new}})$, is almost always cheaper than $C(\tau_{\text{new}}, m)$, the cost of scenario (A). Whenever possible, our basic strategy is therefore to keep $\tau$ constant and use the method described in section 3.6.2 to let $m$ vary up to convergence or to $m_{\max}$. It is only when $m_{\max}$ has been reached that we should consider varying the stepsize to avoid using excessive amount of memory. In the case where $m_{\max}$ is reached, we check if the step is going to be accepted. If it does not satisfy the acceptance criterion, then the procedure given in section 3.6.1 will suggest a new stepsize that is more likely to satisfy the error tolerance. In this case, we use a safety factor $\gamma = 0.6$ to reduce the risk of another rejection. If the tolerance is satisfied, then the stepsize is changed using the default value of $\gamma = 0.9$ to give the possibility to use either a smaller or a larger size in the next step.

Our approach to adaptively propose new stepsize and the dimension is largely inspired by the work of Niesen and Wright [21,22] and it is to be outlined in the next subsection.

**Algorithm 4** Krylov adaptivity.

1: input: A, $(b_0, b_1, \ldots, b_p)$, Tol, m
2: $\gamma_{\max} = 0.6$, $\gamma = 0.9$
3: $\tau = 1$, $\delta = 1.4$
4: **if** happy_breakdown
5: $\quad \omega = 0$
6: $\quad \tau_{\text{new}} = \min(\tau, T_{\text{end}} - T_{\text{now}})$
7: $\quad m_{\text{new}} = m$
8: $\quad$ happy_breakdown = false
9: **else if** $j == m_{\max}$
10: $\quad$ **if** $\omega > \delta$
11: $\qquad m_{\text{new}} = j$
12: $\qquad \tau_{\text{new}} = \tau \left(\frac{\gamma_{\max}}{\omega}\right)^{\frac{1}{q+1}}$
13: $\qquad \tau_{\text{new}} = \min(T_{\text{end}} - T_{\text{now}}, \max(\tau/5, \tau_{\text{new}}))$
14: $\quad$ **else**
15: $\qquad$ Keep $m$ constant and vary $\tau$ as in section 3.6.1
16: $\quad$ **end if**
17: **else**
18: $\qquad$ Vary $m$ as in section 3.6.2 and keep $\tau$ constant
19: **return** $\tau_{\text{new}}$ and $m_{\text{new}}$

### 3.6.1. Variable stepsize

Each interval must be small enough to reduce the norm of the matrix to a level where the problem could be solved with a reasonably large Krylov subspace. On the other hand, it should not be too small because this would necessitate a large number of iterations. The challenge is therefore to find an optimal partition of the interval $[0, \tau]$ into $k$ subintervals $\tau_l$ $(l = 1, \ldots, k)$ without any *a priori* knowledge of the matrix characteristics.

The procedure to determine the stepsize is similar to the adaptive strategy used in many ODE solvers. We assume that the error is approximately $C\tau^{q+1}$, for some constants $q, C \in \mathbb{R}$. The order $q$ is set to $q = \frac{m}{4} - 1$ for the first step as in [21] and if a previously suggested stepsize is rejected, then we use the following estimate to obtain a new order estimate

$$q = \frac{\log(\tau/\tau_{\text{old}})}{\log(\|\epsilon_m\| / \|(\epsilon_m)_{\text{old}}\|)} - 1, \tag{43}$$

where the "old" subscript denotes the rejected estimates.

The suggested stepsize is then given by

$$\tau_{\text{new}} = \tau \left(\frac{\gamma}{\omega}\right)^{\frac{1}{q+1}} \tag{44}$$

where $\omega = T_{\text{end}}\epsilon_m/(\tau \,\text{tol})$ denotes the scaled error and $\gamma$ is a safety factor. Following [8], we chose a value of $\gamma = 0.9$, except when the size of the Krylov subspace $m$ reaches a certain limit $m_{\max}$ in which case we use $\gamma = 0.6$ to avoid step rejections. The suggested stepsize is always clipped to restrict its variation around the current stepsize between $\frac{\tau}{5}$ and $5\tau$ and to make sure it does not exceed the final value of $t_{\text{end}}$.

*3.6.2. Variable dimension of the Krylov subspace*

As mentioned before, a rapid convergence could be obtained if we choose a large value of $m$ when $\|\tau\tilde{A}\|$ is large. Since we do not know the norm of $\tau\tilde{A}$, we can use the *a posteriori* estimate to increase $m$ if the estimated error is too large or to reduce it if the error can be kept small enough [21]. In practical applications, it may also be judicious to select a minimum and maximum size $m_{\min}$ and $m_{\max}$ to cap the memory requirements and to make efficient use of the CPU cache memory. The default values, used in our numerical experiments, are $m_{\min} = 10$ and $m_{\max} = 128$.

The estimate (34) suggests that the error is about $C\kappa^{-m}$ for some constants $C, \kappa \in \mathbb{R}$. For the first suggestion of a step, a default value $\kappa = 2$ is used. When an iteration is rejected, the value of $\kappa$ is estimated as a function of the error estimate from the current and the previously rejected error estimates (again denoted with a "old" subscript) by

$$\kappa = \left(\frac{\omega}{\omega_{\text{old}}}\right)^{1/(m_{\text{old}}-m)} \tag{45}$$

Given this estimate, we obtain the suggested dimension size with

$$m_{\text{new}} = m + \frac{\log(\omega/\gamma)}{\log\kappa} \tag{46}$$

The dimension of the Krylov space is constrained to vary upward or downward in the user specified interval $[m_{\min}, m_{\max}]$. Similarly to what we did for the stepsize, the dimension is clipped to avoid variation larger than 25% below the previous value of $m$ or 33% above.

As an aside, we note that there is a small inconsistency in the methodology to limit the evolution of $m$ and $\tau$ in *phipm*. This algorithm uses maximum and minimum thresholds to avoid large variations of the parameters, but the clipping of the values is done at the end of an iteration. The consequence is that an optimal decision taken on the basis of the cost function may be affected later, leading to suboptimal adaptivity. This is in contrast with KIOPS where the clipping of the suggested values is always done during the adaptive procedure.

*3.7. Avoiding rounding errors*

We close this section with a remark about implementation details that could help to reduce the impact of the rounding errors due to finite precision arithmetic. As shown in Theorem 1, the exact value of $v$ for entries $N+1$ to $N+p$ is given by

$$v(N+1:N+p) = \left[\frac{(T_{\text{now}}+\tau)^{p-1}}{(p-1)!}, \ldots, (T_{\text{now}}+\tau), 1\right]^{\mathsf{T}}. \tag{47}$$

We can also follow the suggestion of Al-Mohy and Higham [6] and normalize the $B$ matrix. To do so, we substitute $B$ by $\nu \cdot B$ and change entries $N+1$ to $N+p$ of the starting vector $v$ to $\mu \cdot v(N+1:N+p)$ before applying Algorithm 2. The normalization constants

$$\nu = 2^{-\lceil \log_2(\|B\|_1)\rceil} \tag{48}$$

$$\mu = 2^{\lceil \log_2(\|B\|_1)\rceil} \tag{49}$$

are defined as a power of 2 to avoid the introduction of rounding errors and has no effect on the solution in exact arithmetic. Our experiments have shown that the normalization of the $B$ matrix leads to slightly faster convergence on most problems.

## 4. Numerical experiments

We evaluate the performance of the KIOPS algorithm in the context of the following fourth- and fifth-order EPIRK exponential methods since these schemes were shown to provide most efficiency in previous studies [33,45,43]. The numerical experiments presented below both validate the performance of the KIOPS algorithm and demonstrate how this algorithm can improve both efficiency and accuracy of the EPIRK and other exponential schemes compared to the adaptive Krylov algorithm *phipm* from [21].

**Table 1**
Coefficients of EPIRK5P1.

$$\begin{bmatrix} \alpha_{11} & & \\ \alpha_{21} & \alpha_{22} & \\ \beta_1 & \beta_2 & \beta_3 \end{bmatrix} = \begin{bmatrix} 0.3512959269505819 & & \\ 0.8440547201165712 & 1.690589160956896 & \\ 1.0 & 1.272712731735689 & 2.271459926542262 \end{bmatrix}$$

$$\begin{bmatrix} g_{11} & & \\ g_{21} & g_{22} & \\ g_{31} & g_{32} & g_{33} \end{bmatrix} = \begin{bmatrix} 0.3512959269505819 & & \\ 0.8440547201165712 & 1.0 & \\ 1.0 & 0.7111109536436687 & 0.6237811195337149 \end{bmatrix}$$

We choose the following four exponential schemes for our numerical experiments:

- EPIRK4s3 [45,43] – stiffly accurate fourth-order integrator:

$$\begin{aligned} U_{n2} &= u_n + \frac{1}{8}\varphi_1\left(\frac{1}{8}h_n J_n\right) h_n f(u_n), \\ U_{n3} &= u_n + \frac{1}{9}\varphi_1\left(\frac{1}{9}h_n J_n\right) h_n f(u_n), \\ u_{n+1} &= u_n + \varphi_1\left(h_n J_n\right) h_n f(u_n) \\ &\quad + \left(1892\varphi_3\left(h_n J_n\right) - 42336\varphi_4\left(h_n J_n\right)\right) h_n r(U_{n2}) \\ &\quad + \left(1458\varphi_3\left(h_n J_n\right) - 34992\varphi_4\left(h_n J_n\right)\right) \\ &\quad \cdot h_n \left(r(U_{n3}) - 2r(U_{n2})\right). \end{aligned} \tag{50}$$

- EPIRK4s3A [33] – stiffly accurate fourth-order integrator:

$$\begin{aligned} U_{n2} &= u_n + \frac{1}{2}\varphi_1\left(\frac{1}{2}h_n J_n\right) h_n f(u_n) \\ U_{n3} &= u_n + \frac{2}{3}\varphi_1\left(\frac{2}{3}h_n J_n\right) h_n f(u_n) \\ u_{n+1} &= u_n + \varphi_1(h_n J_n) h_n f(u_n) \\ &\quad + \left(32\varphi_3(h_n J_n) - 144\varphi_4(h_n J_n)\right) h_n r(U_{n2}) \\ &\quad + \left(-\frac{27}{2}\varphi_3(h_n J_n) + 81\varphi_4(h_n J_n)\right) h_n r(U_{n3}). \end{aligned} \tag{51}$$

- EPIRK5P1 [23] – classical (non-stiffly accurate) fifth-order integrator:

$$\begin{aligned} U_{n2} &= u_n + \alpha_{11}\varphi_1(g_{11} h_n J_n) h_n f(u_n) \\ U_{n3} &= u_n + \alpha_{21}\varphi_1(g_{21} h_n J_n) h_n f(u_n) + \alpha_{22} \\ u_{n+1} &= u_n + \beta_1\varphi_1(g_{31} h_n J_n) h_n f(u_n) + \beta_2\varphi_1(g_{32} h_n J_n) h_n r(U_{n2}) \\ &\quad + \beta_3\varphi_3(g_{33} h_n J_n) h_n(-2r(U_{n2}) + r(U_{n3})) \end{aligned} \tag{52}$$

  with coefficients given in Table 1.
- EXPRB5s3 [49] – stiffly accurate fifth-order method that was originally derived as exponential Rosenbrock integrator and can also be written as an EPIRK scheme [33]:

$$\begin{aligned} U_{n2} &= u_n + \frac{1}{2}\varphi_1(\frac{1}{2}h_n J_n) h_n f(u_n) \\ U_{n3} &= u_n + \frac{9}{10}\varphi_1(\frac{9}{10}h_n J_n) h_n f(u_n) \\ &\quad + \left(\frac{27}{25}\varphi_3(\frac{1}{2}h_n J_n) + \frac{729}{125}\varphi_3(\frac{9}{10}h_n J_n)\right) h_n r(U_{n2}) \\ u_{n+1} &= u_n + \varphi_1(h_n J_n) h_n f(u_n) \\ &\quad + \left(18\varphi_3(h_n J_n) - 60\varphi_4(h_n J_n)\right) h_n r(U_{n2}) \\ &\quad + \left(-\frac{250}{81}\varphi_3(h_n J_n) + \frac{500}{27}\varphi_4(h_n J_n)\right) h_n r(U_{n3}). \end{aligned} \tag{53}$$

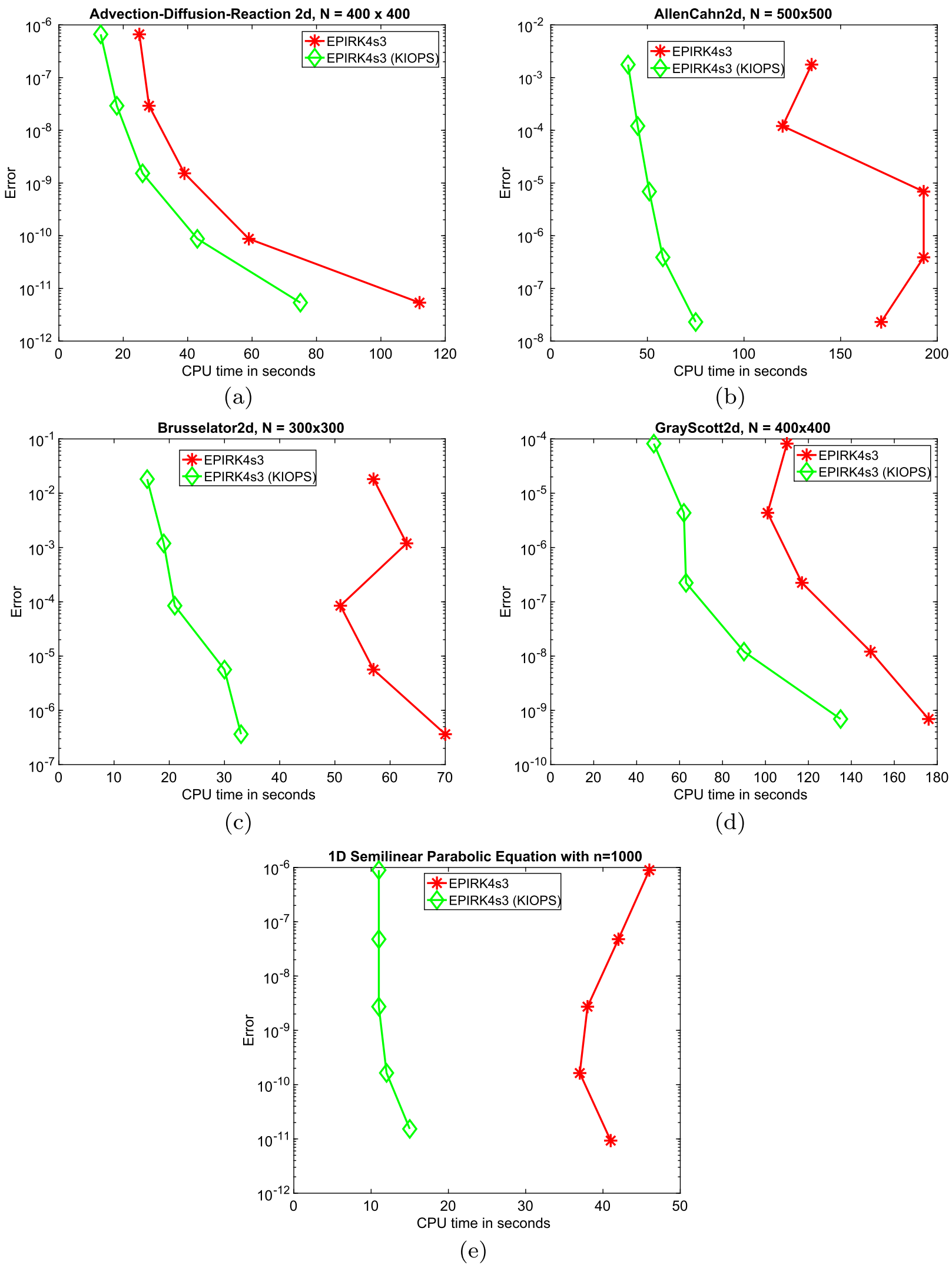

**Fig. 1.** Precision diagram (CPU time vs. error) comparing the adaptive Krylov *phipm* and KIOPS versions of EPIRK4s3 method.

All of these methods are implemented as described in section 2 in a way that groups terms to minimize the number of KIOPS calls needed with each call optimized. A more detailed description of this so-called mixed implementation technique to optimize performance can be found in [33]. With the mixed implementation the fourth-order methods require only two calls to KIOPS or *phipm* routines per time step, while fifth-order methods require only three calls. This is accomplished by simultaneously evaluating the second terms in the right-hand side of equations for $U_{n2}$ and $U_{n3}$ using a single call to KIOPS or *phipm*. For fourth-order methods the second call to KIOPS or *phipm* evaluates the full linear combination of $\varphi$-functions in evaluation of $u_{n+1}$. The fifth-order methods require the third call to KIOPS or *phipm* to approximate the third term in the right-hand side of $U_{n3}$. All schemes are implemented with constant time step integration. The *phipm* algorithm used in our experiments allows for simultaneous evaluation of linear combinations of $\varphi$'s at values $T = [T_1, \ldots, T_{\text{end}}]$. Further details and information on our *phipm* implementation can be found in [21,23,33].

As in previous publications [33], we choose the following test problems routinely used to study the performance of stiff integrators. In all of the problems presented below the $\nabla^2$ term is discretized using the standard second order finite differences.

*Allen–Cahn 2D.* Two-dimensional Allen–Cahn equation [50]:

$$u_t = \alpha \nabla^2 u + u - u^3, \; x, y \in [-1, 1], t \in [0, 1.0]$$

with $\alpha = 0.1$, using no-flow boundary conditions and initial conditions given by $u = 0.1 + 0.1\cos(2\pi x)\cos(2\pi y)$.

*Advection–Diffusion–Reaction (ADR) 2D.* Two-dimensional advection–diffusion–reaction equation [51]:

$$u_t = \epsilon \nabla^2 u - \alpha(u_x + u_y) + \gamma u(u - \tfrac{1}{2})(1 - u), \; x, y \in [0, 1], t \in [0, 0.1],$$

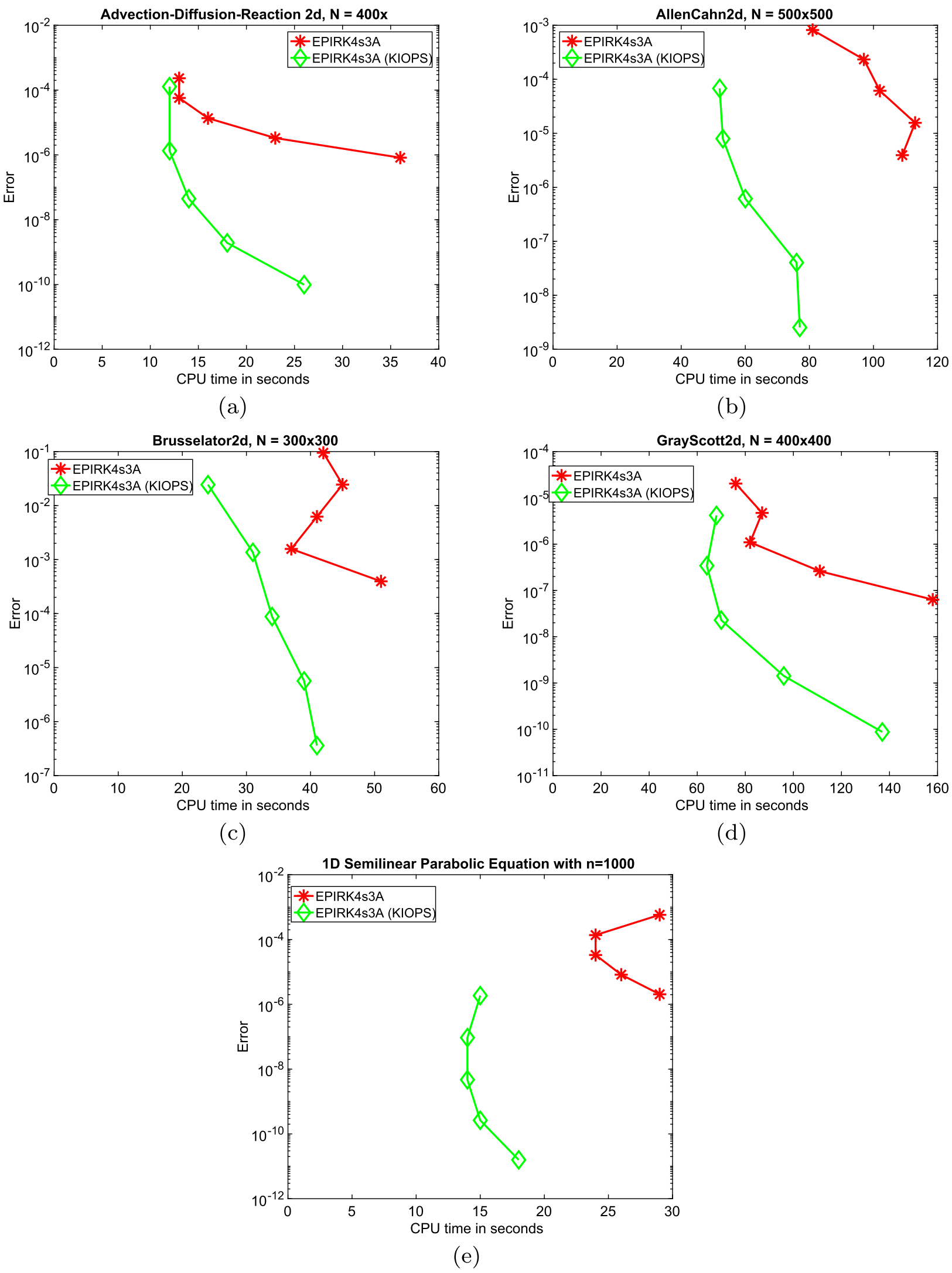

**Fig. 2.** Precision diagram (CPU time vs. error) comparing the adaptive Krylov *phipm* and KIOPS versions of EPIRK4s3A method.

where $\epsilon = 1/100$, $\alpha = -10$, and $\gamma = 100$. Homogeneous Neumann boundary conditions were used and the initial conditions were given by $u = 256(xy(1-x)(1-y))^2 + 0.3$.

*Brusselator 2D.* Two-dimensional Brusselator [52,53]

$$
\begin{aligned}
u_t &= 1 + u^2 v - 4u + \alpha \nabla^2 u, \; x, y \in [0, 1] \\
v_t &= 3u - u^2 v + \alpha \nabla^2 v \\
\alpha &= 0.02
\end{aligned}
$$

with homogeneous Neumann boundary conditions, $t \in [0, 1]$, and initial values

$$
\begin{aligned}
u &= 2 + 0.25y \\
v &= 1 + 0.8x
\end{aligned}
$$

*Gray–Scott 2D.* Two-dimensional Gray–Scott [54] with periodic boundary conditions:

$$
\begin{aligned}
u_t &= d_u \nabla^2 u - uv^2 + a(1-u), \; x, y \in [0, 1] \\
v_t &= d_v \nabla^2 v + uv^2 - (a+b)v,
\end{aligned}
$$

and $d_u = 0.2, d_v = 0.1, a = 0.04, b = 0.06$. Initial conditions given by

$$u = 1 - e^{-150[(x-\frac{1}{2})^2+(y-\frac{1}{2})^2]},$$
$$v = e^{-150[(x-\frac{1}{2})^2+2(y-\frac{1}{2})^2]}.$$

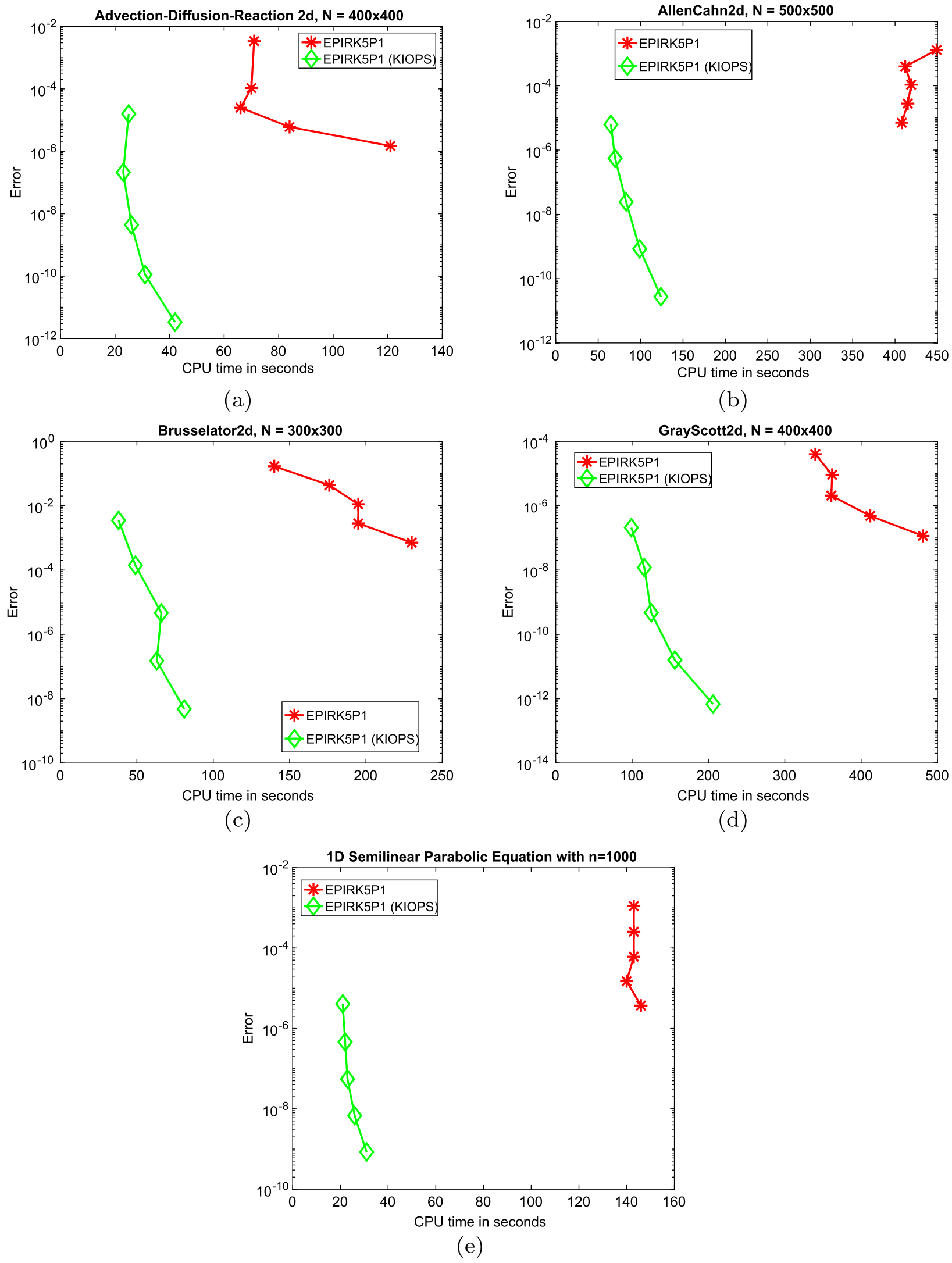

**Fig. 3.** Precision diagram (CPU time vs. error) comparing the adaptive Krylov *phipm* and KIOPS versions of EPIRK5P1 method.

*1D Semilinear Parabolic.* One-dimensional semilinear parabolic problem [55]:

$$\frac{\partial U}{\partial t}(x,t) - \frac{\partial^2 U}{\partial x^2}(x,t) = \int_0^1 U(x,t)dx + \Phi(x,t)$$

with homogeneous Dirichlet boundary conditions and for $x \in [0,1]$ and $t \in [0,1]$. The source function $\Phi$ is chosen such that $U(x,t) = x(1-x)e^t$ is the exact solution.

Figs. 1–4 show precision diagrams (CPU time vs. error) for each of the exponential schemes and all of the test problems comparing the performance of KIOPS and *phipm* versions of the exponential integrators. The following time step sizes and spatial discretizations were used in these simulations

- ADR: $N = 400^2$ with $h = 0.01, 0.005, 0.0025, 0.00125, 6.25 \cdot 10^{-4}$,
- Allen–Cahn: $N = 500^2$ with $h = 0.5, 0.25, 0.1250, 0.0625, 0.03125$,
- Semilinear Parabolic: $N = 1000$ with $h = 0.1, 0.05, 0.0250, 0.0125, 6.25 \cdot 10^{-3}$,
- Gray–Scott: $N = 400^2$ with $h = 0.01, 0.005, 0.0025, 0.00125, 6.25 \cdot 10^{-4}$,
- Brusselator: $N = 300^2$ with $h = 0.25, 0.1250, 0.0625, 0.03125, 0.015625$,

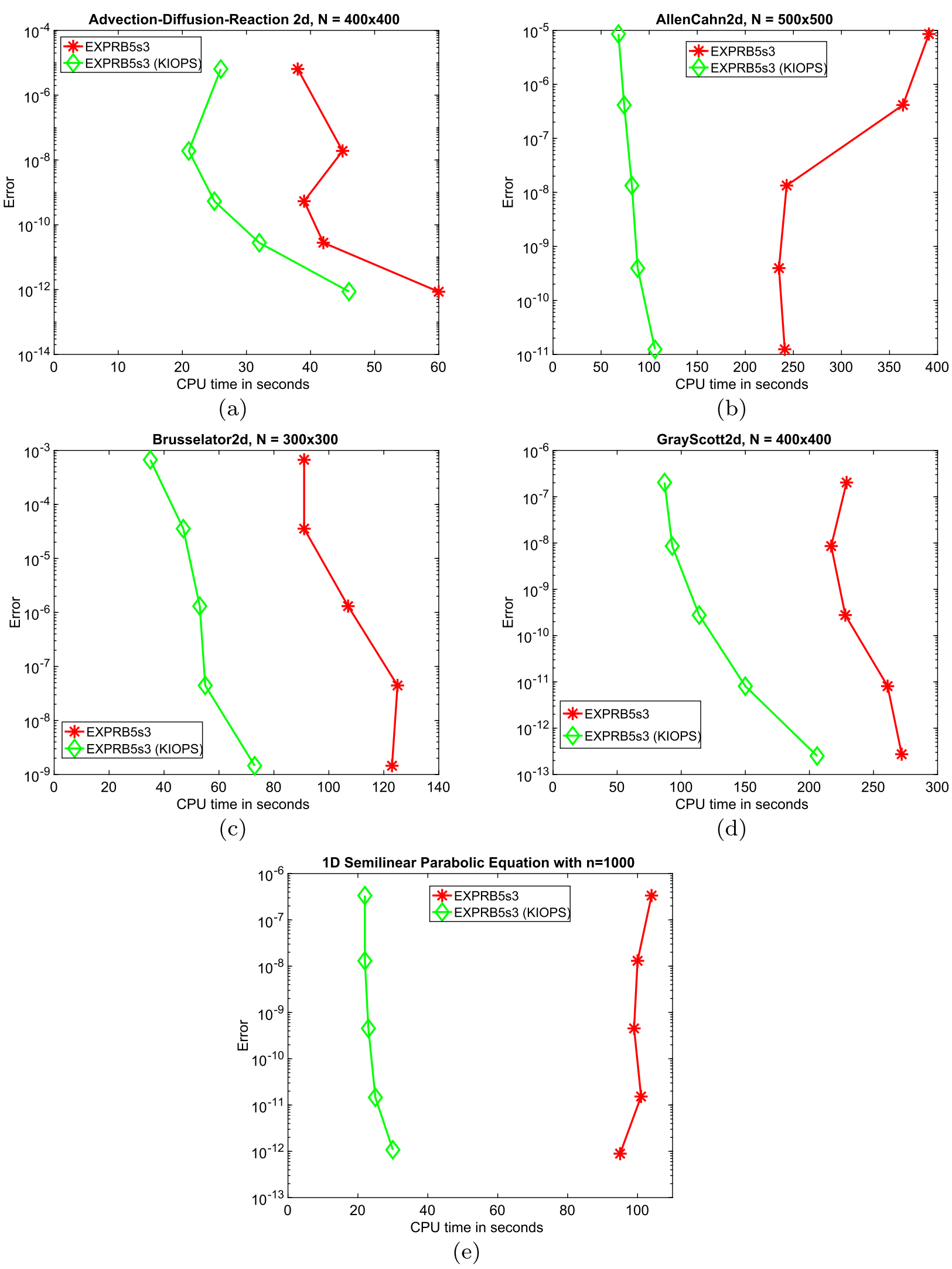

**Fig. 4.** Precision diagram (CPU time vs. error) comparing the adaptive Krylov *phipm* and KIOPS versions of EXPRB5S3 method.

where $N$ and $h$ correspond to the number of spatial discretization nodes and the time step size respectively. The tolerance is set to $10^{-14}$ for KIOPS and *phipm*. The error is defined as the discrete infinity (maximum) norm of the difference between the approximation to the solution and the reference solution computed using MATLAB's *ode15s* integrator with absolute and relative tolerances set to $10^{-14}$.

First, it is important to note that KIOPS outperforms *phipm* in *all* of the simulations delivering both better efficiency and accuracy. For a given tolerance the speedup of simulations can be a factor of 5 or 7 (e.g. Fig. 4(e), Fig. 3(e)). It is interesting to observe that while comparatively EPIRK4s3 is more accurate than EPIRK4s3A and EXPRB5s3 is more accurate than EPIRK5P1 for the *phipm* implementation, KIOPS makes these schemes on par with each other. In other words, a KIOPS EPIRK4s3A implementation is as efficient and accurate as EPIRK4s3, and similar conclusion holds for EXPRB5s3 and EPIRK5P1. It is particularly notable that this phenomenon occurs even between a classically (non-stiffly) accurate EPIRK5P1 and the stiffly accurate method EXPRB5s3. Practice shows that stiff order conditions on exponential integrators can sometimes be unnecessarily strict for some problems [33]. This result raises a question of whether improvements in approximating $\varphi$-vector products can allow for easing strict stiff order conditions and enable derivation of more efficient methods.

To demonstrate that the computational advantage of KIOPS compared to *phipm* is retained as the size of the problem and consequently its stiffness increase we present numerical experiments for each of the problems with different values of the problem size $N$ in Fig. 5. It is expected that as $N$ increases performance of both algorithms will suffer since Krylov subspace projection is at the core of each of these methods and its performance is affected by the more spread out spectrum of the Jacobian matrix. However, graphs in Fig. 5 indicate that the computational advantage of KIOPS compared to *phipm* is retained and even improved as the problem size $N$ increases.

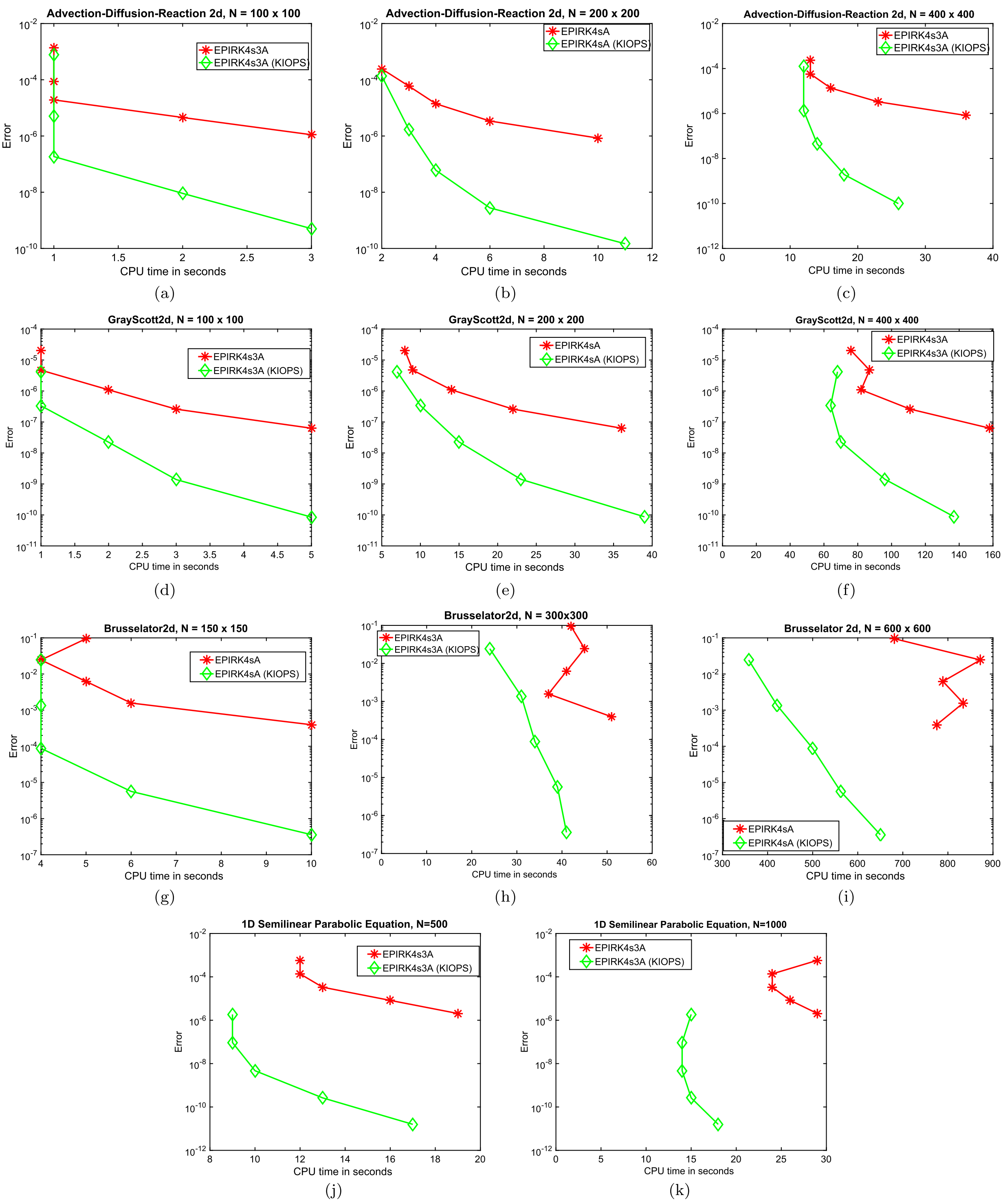

**Fig. 5.** Precision diagram (CPU time vs. error) comparing the adaptive Krylov *phipm* and KIOPS versions of EPIRK4s3A method as the size of each of the problems increases. Advection–Diffusion–Reaction 2D with (a) $N = 100^2$, (b) $N = 200^2$, (c) $N = 400^2$; Gray–Scott 2D with (d) $N = 100^2$, (e) $N = 200^2$; Brusselator 2D with (g) $N = 150^2$, (h) $N = 300^2$, (i) $N = 600^2$; 1D Semilinear Parabolic equation with (j) $N = 500$, (k) $N = 1000$. *Note:* Allen–Cahn problem is omitted since the results for this problem are essentially the same as for other test problems but for clarity we wanted to keep all the graphs in one figure.

## 5. Conclusions

We presented a new KIOPS algorithm for evaluating products of exponential and exponential-like $\varphi$-functions of large stiff matrices and vectors. To date, *phipm* was considered to be the most efficient algorithm for problems where no additional information about the spectrum or norm of the matrix is available. Our results demonstrated that the new KIOPS algorithm outperforms *phipm*. The efficiency of the proposed algorithm is attributable to a combination of the incomplete orthogonalization procedure, a better adaptivity procedure and a heuristic strategy to determine the initial size of the Krylov space using information from previous time substeps. The new algorithm offers not only better computational efficiency but new pathways to further improvements in making exponential and exponential-type integrators more computationally appealing. In particular, the adaptivity algorithms within the KIOPS method can be improved further if better error and cost estimates are derived. We plan to pursue this line of research in our future work. We note that alternative techniques, like restarted Krylov subspace or block Krylov subspace, could be combined with the adaptive method described in this paper. We intend to study them in future work. We also plan to investigate effective ways to improve the performance of the algorithm on parallel architectures.

## 6. Code availability

The EPIC package implements the exponential integrators used in our numerical experiments. It is available from http://faculty.ucmerced.edu/mtokman/#software.

## Acknowledgements

The origin of this study goes back to fruitful discussions with Michel Desgagné and Martin Charron. Stimulating discussions with Janusz Pudykiewicz and Michel Valin are also acknowledged. Reviews by Christopher Subich, Rabah Aider, and Nathalie Berger greatly improved the original manuscript. This work was in part supported by a grant no. 1115978 from the National Science Foundation, Computational Mathematics Program.

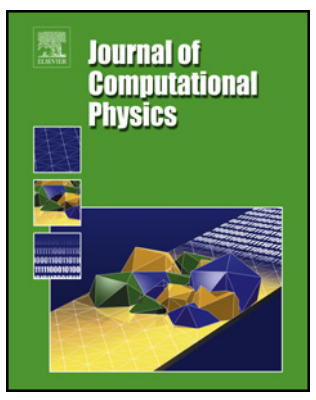

Corrigendum

# Corrigendum to “KIOPS: A fast adaptive Krylov subspace solver for exponential integrators” [J. Comput. Phys. 372 (2018) 236–255]

Stéphane Gaudreault [a,*], Greg Rainwater [b], Mayya Tokman [b]

[a] *Recherche en prévision numérique atmosphérique, Environnement et Changement climatique Canada, 2121 Route Transcanadienne, Dorval, Québec, H9P 1J3, Canada*

[b] *School of Natural Sciences, University of California, 5200 N. Lake Road, Merced, CA 95343, United States*

1. Minor errors were present in Eqs. (18) and (19) of [1]:
   “and the columns of $E_{12}$ are given by the Theorem 2.1 of Al-Mohy and Higham [6], with $l = 0$ and noting our reverse ordering of the columns in $B$, as

$$E_{12}(1:N, N+j) = \sum_{k=1}^{j} \tau^k \varphi_k(\tau A) b_{p-j+k}\,, \quad \text{for } j = 1 \text{ to } p. \tag{18}$$

   Inserting Eqs. (17) and (18) into Eq. (16), we obtain the following matrix”

$$e^{\tau \tilde{A}} = \begin{bmatrix} e^{\tau A} & \tau\varphi_1(\tau A) b_p & \tau\varphi_1(\tau A) b_{p-1} + \tau^2 \varphi_2(\tau A) b_p & \dots & \sum_{k=1}^{p} \tau^k \varphi_k(\tau A) b_k \\ 0 & 1 & \frac{\tau}{1!} & \dots & \frac{\tau^{p-1}}{(p-1)!} \\ \vdots & & \ddots & \ddots & \vdots \\ & & & 1 & \frac{\tau}{1!} \\ 0 & & & 0 & 1 \end{bmatrix} \tag{19}$$

2. A typographical error was found in Algorithm 1 where line 13 should be “$w(N+1:N+p, l)$”.

* Corresponding author.
*E-mail addresses:* stephane.gaudreault2@canada.ca (S. Gaudreault), grainwater@ucmerced.edu (G. Rainwater), mtokman@ucmerced.edu (M. Tokman).

None of the corrections presented here affect the conclusions of the original article, nor the results of the numerical experiments.